\documentclass[11pt]{amsart}

\usepackage[margin=1.23in]{geometry}

\usepackage[T1]{fontenc}
\usepackage{lmodern}
\usepackage[utf8]{inputenc}
\usepackage{amsmath,amssymb,amsthm,mathtools,mathrsfs}
\usepackage{enumitem}
\usepackage{xcolor}
\usepackage{microtype}
\usepackage[colorlinks=true,linkcolor=blue,citecolor=blue,urlcolor=blue]{hyperref}

\setlist[enumerate]{label=(\roman*),leftmargin=2.2em}
\setlist[itemize]{leftmargin=2.2em}
\newcommand{\gen}[2]{\left\langle #1 \right\rangle_{#2}}

\newtheorem{theorem}{Theorem}[section]
\newtheorem{lemma}[theorem]{Lemma}
\newtheorem{proposition}[theorem]{Proposition}
\newtheorem{corollary}[theorem]{Corollary}

\theoremstyle{definition}
\newtheorem{definition}[theorem]{Definition}
\newtheorem{example}[theorem]{Example}
\newtheorem{notation}[theorem]{Notation}
\theoremstyle{remark}
\newtheorem{remark}[theorem]{Remark}

\newcommand{\SU}{\operatorname{SU}}
\newcommand{\acl}{\operatorname{acl}}
\newcommand{\dcl}{\operatorname{dcl}}
\newcommand{\tcl}{\operatorname{tcl}} 
\newcommand{\tp}{\operatorname{tp}}
\newcommand{\Th}{\operatorname{Th}}
\newcommand{\Aut}{\operatorname{Aut}}
\newcommand{\Def}{\operatorname{Def}}
\newcommand{\Comp}{\operatorname{Comp}}
\newcommand{\meas}{\operatorname{meas}}
\newcommand{\dimn}{\operatorname{dim}}

\newcommand{\MS}{\mathrm{MS}}

\newcommand{\N}{\mathbb{N}}

\newcommand{\R}{\mathbb{R}}

\newcommand{\dsum}{\mathbin{\dot\cup}}
\newcommand{\concat}{\mathbin{{}^{\smallfrown}}}
\newcommand{\eps}{\varepsilon}
\newcommand{\emptyseq}{\langle\rangle}
\newcommand{\loc}{\operatorname{loc}}

\title[MS-measurability via coordinatization]{MS-measurability via coordinatization}
\author{Mostafa Mirabi}
\address{Department of Mathematics and Computer Science, Wesleyan University, Middletown, CT 06459, USA, and The Taft School, Watertown, CT 06795, USA}

\email{mmirabi@wesleyan.edu}
\urladdr{https://sites.google.com/site/mostafamirabi/}


\begin{document}
\raggedbottom

\begin{abstract}
We study $\MS$-measurability of structures obtained as expansions of finite-height
trees. We show that an admissible nil-interaction coordinatized structure is
$\MS$-measurable whenever its canonical component theories are $\MS$-measurable
and their universes have dimension one. The dimension of a complete type is the
sum of its local component dimensions, and its measure is the product of the
corresponding local measures.

We also show that, in the countable $\aleph_0$-categorical setting, an admissible
nil-interaction structure is interdefinable with its canonical tree product. These products
preserve supersimplicity of finite $\SU$-rank and one-basedness, with global
$\SU$-rank given by a sum of local ranks. Finite homogeneous component
approximations give finite tree envelopes, yielding preservation of smooth
approximability. Enriched finite tree products of multidimensional asymptotic
classes, respectively exact classes, again form multidimensional asymptotic
classes, respectively exact classes. For the corresponding reducts, the
conclusion is weak unless the coordinate enrichment can be recovered uniformly.
\end{abstract}

\maketitle

\section{Introduction}

The notion of $\MS$-measurability was introduced by Macpherson and Steinhorn
in \cite{MacphersonSteinhorn2008}. An $\MS$-measurable structure is an infinite
structure whose definable sets are equipped with an $\N$-valued dimension and
an $\R$-valued measure satisfying finiteness, definability, and Fubini axioms.
We recall these axioms in Section~\ref{sec:omega-cat-criterion}.
The motivating examples come from finite fields. The theorem of Chatzidakis,
van den Dries and Macintyre gives uniform estimates for the sizes of definable
sets over finite fields: the size is asymptotic to a rational multiple of a
power of the field size, and the dimension and measure vary definably with
the parameters \cite{CDM1992}.

Macpherson and Steinhorn also introduced one-dimensional asymptotic classes
of finite structures. Roughly speaking, an asymptotic class is a class in
which the sizes of definable sets satisfy uniform estimates of this kind.
Elwes developed higher-dimensional asymptotic classes and showed that their
infinite ultraproducts are measurable; see
\cite{ElwesThesis,Elwes2007,ElwesMacphersonSurvey}. More recently, the theory
of multidimensional asymptotic and exact classes has allowed several
independent size parameters and generalized measures with values in a
semiring \cite{AnscombeMacphersonSteinhornWolf2024}. A related approach to
simplicity through pseudofinite dimension is studied in
\cite{GarciaMacphersonSteinhorn2015}.

Coordinatization constructs a structure from a skeleton and local components,
with prescribed relations between them. Here the skeleton is a finite-height
tree. Cherlin and Hrushovski developed Lie coordinatization and finite
envelopes in their study of smoothly approximable structures
\cite{CherlinHrushovski2003}. Related methods have been used to construct
pseudofinite structures and asymptotic classes \cite{HillPreprint,Hill2017}.
The pure trees used here are described by finite tree plans, whose asymptotic
properties are studied in \cite{Mirabi2025Trees}.

We ask when measurability of the components implies measurability of the
whole structure. Relations between different components are the main
obstacle. For example, place the two parts of a generic bipartite graph above
a common root. An edge to a parameter in the other part changes the type of
a point without changing its local type in its own component. A construction
using only local dimensions and measures therefore requires a restriction on
these interactions.

In this paper, we study a strong notion of coordinatization called
\emph{nil-interaction}. Our main result is the following.

\begin{quote}
Let $M$ be an admissible nil-interaction expansion of a finite tree plan
$\Gamma$. Suppose that, for each infinity-node $\sigma$ of $\Gamma$, the
canonical component theory $T_\sigma$ is $\MS$-measurable, with its universe
normalized to have dimension one. Then $M$ is $\MS$-measurable. If $C$ is
finite and tree-closed and $a$ is a finite tuple, then
$$
\begin{aligned}
   \dimn(a/C)&=\sum_K\dimn_{T_{\sigma(K)}}(\bar a_K/C_K),\\
   \meas(a/C)&=\prod_K\meas_{T_{\sigma(K)}}(\bar a_K/C_K),
\end{aligned}
$$
where $K$ ranges over the active components in $\tcl(Ca)$ over $C$, and
$\bar a_K$ enumerates the new elements in $K$.
\end{quote}

A \emph{tree plan} is a finite rooted tree whose nodes are labelled by $1$
or $\infty$. A $1$-node gives a unique successor of the corresponding
plan-type; an $\infty$-node gives a copy of a prescribed component structure
above each appropriate predecessor. The resulting \emph{canonical
$\Gamma$-tree product} is presented in a one-sorted language containing the
tree structure and the local component relations.

Nil-interaction says that, over a finite tree-closed set $C$, the type of
a new point $e$ is determined by its tree position and its local type over
$C\cap\Comp(e)$. We add $e$ together with its forced $1$-successors, so each
stage remains tree-closed. A common canonical language allows comparison
of local types in different copies of the same component.

Adding these blocks one at a time gives a normal form: the type of a finite
tuple is determined by its pure tree type and finitely many local component
types. It follows that every countable admissible nil-interaction structure
is interdefinable with its canonical tree product. The theorem thus concerns
this form of coordinatization; it does not classify all one-based or all
measurable $\aleph_0$-categorical structures.

Finiteness of the plan ensures that the closure of a finite tuple meets
only finitely many components, and tree-closed bases contain all forced
successors. The assumption $\acl=\tcl$ excludes further algebraic points.
By $\aleph_0$-categoricity, complete types over finite parameter sets are
isolated. We define dimension and measure on these types and then extend
them to all definable sets.

Pure-set components give the pure finite-height trees, while random-graph
components give unstable $\MS$-measurable examples. Vector spaces over
finite fields are measurable but have linear span as algebraic closure, so
they fall outside the present closure assumption. They and measurable
modules \cite{Kestner2014} would require a closure combining tree closure
with algebraic closure inside the components.

Canonical tree products also preserve supersimplicity of finite $\SU$-rank
and one-basedness. The rank is computed by a local sum:
$$
   \SU_M(a/C)=
   \sum_{K\in\mathcal K(\tcl(Ca)/C)}
   \SU_{T_{\sigma(K)}}(\bar a_K/C_K).
$$
These questions are related to the study of trivial dependence and finite
satisfiability in $\aleph_0$-categorical structures; see
\cite{Djordjevic2006}.

Finite fully homogeneous component approximations combine into finite fully
homogeneous substructures of the product, which we call \emph{tree envelopes}.
This proves preservation of smooth approximability. The construction uses
independent component approximations; its relation to Cherlin--Hrushovski
$\mu$-envelopes is discussed in Section~\ref{sec:rank-envelopes}.

For finite tree products, independent component sizes lead to the
multidimensional asymptotic setting. Enriched products preserve
multidimensional asymptotic and exact classes. The one-sorted reducts satisfy
the corresponding weak conclusions, with full definability when the
enrichment is uniformly recoverable. This result is independent of the
infinite $\MS$-measurability theorem.

Sections~\ref{sec:trees} and~\ref{sec:nic} develop the tree construction,
nil-interaction, and the normal form. Sections~\ref{sec:omega-cat-criterion}
and~\ref{sec:lifting} give the measurability criterion and prove the main
theorem. Section~\ref{sec:rank-envelopes} treats rank, one-basedness, and
tree envelopes; Section~\ref{sec:mac-trees} treats finite products. The final
section gives directions for further work.

\section{Tree plans and tree products}\label{sec:trees}

\subsection{Finite tree plans}

We work in the tree language $L_t=\{\leq,\eps,\wedge,\operatorname{pred}\},$
where $\leq$ is a binary relation symbol, $\eps$ is a constant symbol,
$\wedge$ is a binary function symbol, and $\operatorname{pred}$ is a unary
function symbol. Throughout this paper, all tree plans are finite.

\begin{definition}
\begin{enumerate}
\item A \emph{finite rooted tree} is a finite $L_t$-structure $A$ such that $\leq^A$ is a partial order, $\eps^A$ is the unique minimum element, $a\wedge^A b$ is the greatest element below both $a$ and $b$, $\operatorname{pred}^A(\eps^A)=\eps^A$, and if $a\ne\eps^A$, then $\operatorname{pred}^A(a)$ is the greatest element strictly below $a$.
\item A \emph{tree plan} is a pair $\Gamma=(\Gamma,\lambda)$, where $\Gamma\subseteq\omega^{<\omega}$ is a finite rooted subtree containing $\emptyseq$ and
$$
\lambda:\Gamma\longrightarrow\{1,\infty\}
$$
satisfies $\lambda(\emptyseq)=1$.
\end{enumerate}
\end{definition}

We write $I(\Gamma)=\{\sigma\in\Gamma:\lambda(\sigma)=\infty\}$
for the set of infinity-nodes of the plan.  If $\sigma\in\Gamma$, its height is denoted by $h(\sigma)$.

\begin{definition}
Let $\Gamma=(\Gamma,\lambda)$ be a tree plan and let $X$ be a nonempty set.  The tree $\Gamma(X)$ consists of the root $\emptyseq$ together with all finite
sequences
$\big\langle (i_0,t_0),\ldots,(i_{m-1},t_{m-1})\big\rangle$
 where $m>0$, $\langle i_0,\ldots,i_{m-1}\rangle\in\Gamma$,
and, for each $k<m$,
$$
t_k=\star \quad\text{if }
   \lambda(\langle i_0,\ldots,i_k\rangle)=1,
$$
whereas
$$
t_k\in X \quad\text{if }
   \lambda(\langle i_0,\ldots,i_k\rangle)=\infty.
$$
The map
$$
\pi_X:\Gamma(X)\longrightarrow \Gamma
$$
sends $\emptyseq$ to $\emptyseq$. For a non-root node, it sends
$\big\langle (i_0,t_0),\ldots,(i_{m-1},t_{m-1})\big\rangle$ to
$\langle i_0,\ldots,i_{m-1}\rangle$.  The tree language is interpreted by
initial-segment order, root, meet, and predecessor.  For $\sigma\in\Gamma$ let 
$P_\sigma(\Gamma(X))=\pi_X^{-1}(\sigma).$
Thus $P_\sigma$ denotes the set of nodes whose image under $\pi_X$ is
$\sigma$. We call $\sigma$ their plan-type. We regard $\Gamma(X)$ as an
$L_\Gamma$-structure, where
$L_\Gamma=L_t\cup\{P_\sigma:\sigma\in\Gamma\}.$
\end{definition}

The following figure shows a tree plan $(\Gamma,\lambda)$ and the
corresponding tree $\Gamma(X)$ for $X=\{a,b,c\}$.

\begin{center}
\begin{picture}(350,90)

        \qbezier[50](30,5)(30,25)(30,35)
        \put(30,5){\circle*{3}}
        \put(30,35){\circle*{3}}
        \qbezier[50](30,35)(20,50)(10,65)
        \put(10,65){\circle*{3}}
        \qbezier[50](30,35)(40,50)(50,65)
        \put(50,65){{\color{red}\circle*{5}}}

        \put(34,3){\footnotesize{$\gen{}{}$}}
   			\put(22,3){\footnotesize{{\color{blue}$1$}}}

        \put(34,30){\footnotesize{$\gen{0}{}$}}
   			\put(22,30){\footnotesize{{\color{blue}$1$}}}

        \put(10,70){\footnotesize{$\gen{0,0}{}$}}
   			\put(2,60){\footnotesize{{\color{blue}$1$}}}

        \put(53,70){\footnotesize{$\gen{0,1}{}$}}
   			\put(37,63){\footnotesize{{\color{red}$\infty$}}}

        \put(20,-17){\footnotesize{$(\Gamma, \lambda)$  }}

        \put(160,5){\circle*{3}}
                   \put(100,65){\circle*{3}}
                      \qbezier[100](160,35)(130,50)(100,65)

       \put(160,35){\circle*{3}}
       \qbezier[50](160,5)(160,20)(160,35)

        \put(200,65){{\color{red}\circle*{5}}}
        \put(260,65){{\color{red}\circle*{5}}}
        \put(320,65){{\color{red}\circle*{5}}}

       \qbezier[100](160,35)(180,50)(200,65)
       \qbezier[150](160,35)(210,50)(260,65)
       \qbezier[220](160,35)(240,50)(320,65)

        \put(164,3){\footnotesize{$\gen{}{}$}}
        \put(125,30){\footnotesize{$\gen{(0,\star)}{}$}}
        \put(84,70){\footnotesize{$\gen{(0,\star), (0,\star)}{}$}}
        \put(168,70){\footnotesize{$\gen{(0,\star),(1,a)}{}$}}
   \put(240,70){\footnotesize{$\gen{(0,\star),(1,b)}{}$}}
   \put(310,70){\footnotesize{$\gen{(0,\star),(1,c)}{}$}}
        \put(140,-17){\footnotesize{$\Gamma(\{a,b,c\})$}}

        \end{picture}

\end{center}

\bigskip
\bigskip

For the infinite structures considered below, we use the same countably
infinite set $S$ at every infinity-node. In
Section~\ref{sec:mac-trees}, we allow the finite component sets to vary
with the plan node.

\subsection{Tree closure and components}

\begin{definition}
Let $B\subseteq\Gamma(X)$.  Define $\tcl_0(B)$ by
$$
\tcl_0(B)=\{\eps\}\cup\{a\in\Gamma(X):(\exists b\in B)\ a\leq b\}.
$$
Inductively, let 
$$
\tcl_{n+1}(B)=\tcl_n(B)\cup
   \{a\in\Gamma(X):\lambda(\pi(a))=1\text{ and }\operatorname{pred}(a)\in\tcl_n(B)\}.
$$
Finally,
$$
\tcl(B)=\bigcup_{n<\omega}\tcl_n(B).
$$
A set $B$ is \emph{tree-closed} if $B=\tcl(B)$.  If $a\in\Gamma(X)$ and $B\subseteq\Gamma(X)$, define
$[a\wedge B;X]=\max\{e\in\tcl(B): e\leq a\}.$
When $X$ is fixed, we write $[a\wedge B]$.
\end{definition}

Since $\Gamma$ is finite, $\tcl(B)$ is finite whenever $B$ is finite.

\begin{notation}
For $a\in\Gamma(X)$ define its \emph{component}
$$
\Comp(a)=\{a'\in\Gamma(X):\operatorname{pred}(a')=\operatorname{pred}(a)\text{ and }\pi(a')=\pi(a)\}.
$$
If $\pi(a)\notin I(\Gamma)$, then $\Comp(a)=\{a\}$.  If $\pi(a)\in I(\Gamma)$, then $\Comp(a)$ is a copy of the underlying set $X$ over the predecessor of $a$.
\end{notation}

Adding an element also adds its forced $1$-successors under tree closure.
The normal-form proof will proceed through such tree-closed extensions.

\subsection{The canonical tree product}

We expand $\Gamma(S)$ by placing a copy of a fixed structure $N_\sigma$
in each component of plan-type $\sigma$. Local relations include the
predecessor as an extra coordinate, giving a one-sorted presentation.

\begin{definition}
Let $\Gamma=(\Gamma,\lambda)$ be a tree plan, and let $S$ be a countably infinite set.
For each $\sigma\in I(\Gamma)$ let $N_\sigma$ be a countably infinite relational
structure with universe $S$, in a relational language $L_\sigma$.  Assume that the
languages $L_\sigma$ are pairwise disjoint.
The \emph{canonical $\Gamma$-tree product}
$
\prod_{\sigma\in I(\Gamma)}^{\!*}N_\sigma
$
is the expansion $M$ of $\Gamma(S)$ in the language
$$
L^*
   =
   L_\Gamma\cup
   \{\widehat R:R\in L_\sigma\text{ for some }\sigma\in I(\Gamma)\},
$$
where, if $R\in L_\sigma$ is $n$-ary, then $\widehat R$ is an $(n+1)$-ary
relation symbol.  Suppose $\sigma=\tau\concat\langle k\rangle$.  We interpret
$\widehat R$ by
$$
\begin{aligned}
\widehat R^M
=
\bigg\{&(b,b\concat\langle(k,s_0)\rangle,\ldots,
      b\concat\langle(k,s_{n-1})\rangle): \\
&\pi(b)=\tau
\text{ and }
N_\sigma\models R(s_0,\ldots,s_{n-1})\bigg\}.
\end{aligned}
$$
Thus, over each node $b$ of plan-type $\tau$, the set $\left\{b\concat\langle(k,s)\rangle:s\in S\right\}$
is a copy of the component structure $N_\sigma$.  Distinct components have no
relations between them except those already coming from the pure tree structure.
\end{definition}

Notice that the first coordinate of $\widehat R$ specifies the predecessor.
In particular, the local component structures are uniformly definable over
their predecessors.

\begin{example}
Suppose that every $N_\sigma$ is the countable random graph. Then every
infinite component is a random graph, and there are no additional relations
between different components. Thus the construction allows unstable
components while keeping their interaction nil.
\end{example}

\begin{example}
Let $\Gamma$ consist of a root and two infinity successors, and expand the
two infinite parts by a generic bipartite edge relation. Let $x$ belong to
the left part and let $b$ belong to the right part. The formulas $E(x,b)$
and $\neg E(x,b)$ distinguish points with the same pure tree type and the
same local type in the left component. Hence this presentation of the
generic bipartite graph does not satisfy nil-interaction. A notion of
coordinatization which includes this example must allow relations between
different components.
\end{example}

\section{Nil-interaction coordinatization}\label{sec:nic}

We define nil-interaction, add the assumptions of $\aleph_0$-categoricity
and algebraic triviality, and prove the normal form used in the measuring
construction.

\subsection{Canonical component languages}

Let $M$ be an expansion of $\Gamma(S)$ in a language extending $L_\Gamma$.
If $e\in M$ and $\pi(e)=\sigma\in I(\Gamma)$, let 
$b=\operatorname{pred}(e)$.  The set $\Comp(e)$ is definable over $b$.

Components of the same plan-type will be compared in the following common
language.

\begin{definition}
Let $\sigma\in I(\Gamma)$.  The \emph{canonical component language}
$L_\sigma^{\mathrm{can}}$ contains, for each $n\geq 1$ and each complete
$\emptyset$-type $p(z,x_0,\ldots,x_{n-1})$
realized in $M$ by tuples satisfying
$$
\pi(x_i)=\sigma,\qquad z=\operatorname{pred}(x_i),
   \qquad x_i\in\Comp(x_0)
   \quad\text{for all }i<n,
$$
an $n$-ary relation symbol $R_p$.

If $e\in M$, $\pi(e)=\sigma$, and $b=\operatorname{pred}(e)$, the
\emph{canonical component structure} $M[e]$ is the
$L_\sigma^{\mathrm{can}}$-structure with universe $\Comp(e)$ in which
$$
R_p^{M[e]}(a_0,\ldots,a_{n-1})
   \quad\Longleftrightarrow\quad
   (b,a_0,\ldots,a_{n-1})\models p.
$$
\end{definition}

Note that $L_\sigma^{\mathrm{can}}$ depends on $\sigma$ and $M$, and is the
same for every choice of $e$ of plan-type $\sigma$. If $M$ is
$\aleph_0$-categorical, the types in the definition are isolated. Therefore
the canonical component relations are definable in $M$.

For a tree product, this canonical structure is interdefinable with the
original component structure. The common language allows us to state the
comparison of local types precisely.

\begin{definition}
Let $C\subseteq M$ be finite and tree-closed.  Suppose $e,e'\in M$ and
$\pi(e)=\pi(e')=\sigma\in I(\Gamma).$
We write $e\equiv_C^{\loc}e'$
if the following local comparison holds.
\begin{enumerate}
\item If $\Comp(e)=\Comp(e')$, then
$\tp_{M[e]}(e/C\cap\Comp(e)) = \tp_{M[e]}(e'/C\cap\Comp(e)).$

\item If $\Comp(e)\ne\Comp(e')$, then $C\cap\Comp(e)=C\cap\Comp(e')=\emptyset,$
and $e$ and $e'$ have the same complete $1$-type in the common language
$L_\sigma^{\mathrm{can}}$.
\end{enumerate}
\end{definition}

For a canonical tree product whose component theory has a unique complete
$1$-type, the second clause holds automatically.

\subsection{The local nil-interaction condition}

For a tree-closed set $C$ and an element $e$, write $B(C,e)=\tcl(Ce)$
for the tree-closed set obtained by adding $e$ to $C$. If
$\operatorname{pred}(e)\in C$, $\pi(e)\in I(\Gamma)$, and $e\notin C$, then
$B(C,e)\setminus C$ consists of $e$ and its forced $1$-successors. We call
this a block and use these blocks in the induction below.

\begin{definition}
Let $M$ be an expansion of $\Gamma(S)$.  We say that $M$ satisfies
\emph{nil-interaction over $\Gamma$} if the following two conditions hold.

\begin{enumerate}[label=\textup{NI\arabic*.},leftmargin=4.0em]
\item \textup{(Uniform components)}
If $e,e'\in M$ have the same pure tree $1$-type, that is, $\tp_{\Gamma(S)}(e)=\tp_{\Gamma(S)}(e'),$
then
$\tp_M(e)=\tp_M(e').$

\item \textup{(Local one-block extension)}
Let $C\subseteq M$ be finite and tree-closed, and let $e,e'\in M\setminus C$
satisfy
$$
\pi(e)=\pi(e')\in I(\Gamma),
   \qquad
   \operatorname{pred}(e),\operatorname{pred}(e')\in C.
$$
Suppose there is an isomorphism of pure $L_\Gamma$-trees
$$
\iota:B(C,e)\longrightarrow B(C,e')
$$
fixing $C$ pointwise and sending $e$ to $e'$.  Suppose also that
$e\equiv_C^{\loc}e'.$
Then $\tp_M(B(C,e)/C)=\tp_M(B(C,e')/C),$
where the two blocks are enumerated in corresponding order under $\iota$.
\end{enumerate}
\end{definition}

NI2 applies to the tree-closed block extensions used in the next subsection.

\begin{definition}
An expansion $M$ of $\Gamma(S)$ is an \emph{admissible nil-interaction
coordinatized structure}, or \emph{admissible NIC structure}, if
\begin{enumerate}
\item $M$ satisfies nil-interaction over $\Gamma$;
\item $M$ is countably infinite and $\aleph_0$-categorical;
\item for every finite $A\subseteq M$,
$\acl_M(A)=\tcl(A).$
\end{enumerate}
\end{definition}

We do not assume quantifier elimination in this definition. For the lifting
theorem, we use the normal form below and the equality $\acl=\tcl$.
We will later prove quantifier elimination for canonical tree products
under the stated assumptions on their components.

\subsection{Finite-type normal form}

Let $C\subseteq M$ be finite and tree-closed, and let $A\supseteq C$ be finite
and tree-closed.  A component $K$ is \emph{active in $A/C$} if
$$
K=\Comp(e)
   \quad\text{for some }e\in A\setminus C\text{ with }\pi(e)\in I(\Gamma).
$$
For such a component let
$
A_K=A\cap K,
   \text{ and }
   C_K=C\cap K.
$
If $A'$ is another finite tree-closed set and $K'$ is a component, we similarly
write
$
A'_{K'}=A'\cap K',
   \text { and }
   C_{K'}=C\cap K'.
$

\begin{lemma}\label{lem:block-enumeration}
Let $C\subseteq A$ be finite tree-closed subsets of $\Gamma(S)$.  There are finite
tree-closed sets $C=A_0\subseteq A_1\subseteq\cdots\subseteq A_m=A$
and elements $e_i\in A_i\setminus A_{i-1}$ such that
$$
\pi(e_i)\in I(\Gamma),
   \qquad
   \operatorname{pred}(e_i)\in A_{i-1},
   \qquad
   A_i=\tcl(A_{i-1}e_i)
$$
for each $1\leq i\leq m$.
\end{lemma}

\begin{proof}
We construct the sets inductively. Let  $A_0=C$, and suppose that
$A_i\neq A$. Choose a minimal element $u\in A\setminus A_i$ with respect
to $\leq$. Since $A$ is tree-closed, $\operatorname{pred}(u)\in A$.
Minimality of $u$ gives $\operatorname{pred}(u)\in A_i$.

We claim that $\lambda(\pi(u))=\infty$. Otherwise,
$\lambda(\pi(u))=1$, so tree closure of $A_i$ would imply $u\in A_i$.
This is a contradiction. Hence $\pi(u)\in I(\Gamma)$.

Set $e_{i+1}=u$ and
$A_{i+1}=\tcl(A_i e_{i+1}).$
Since $A$ is tree-closed, $A_{i+1}\subseteq A$. At each step we add
at least one element of $A$. Thus, after finitely many steps, we obtain $A$.
\end{proof}

\begin{theorem}\label{thm:normal-form}
Let $M$ be an admissible NIC structure.  Let $C\subseteq M$ be finite and
tree-closed, and let $A,A'$ be finite tree-closed subsets containing $C$.  Suppose
$
f:A\longrightarrow A'
$ 
is an isomorphism of pure $L_\Gamma$-trees fixing $C$ pointwise.

For each component $K$ active in $A/C$, let $f(K)$ denote the component
containing $f(e)$ for any $e\in A_K\setminus C_K$.  Suppose that, for every
component $K$ active in $A/C$, the restriction
$$
f\restriction A_K:A_K\longrightarrow A'_{f(K)}
$$
is a partial elementary map between the corresponding canonical component
structures, viewed in the common language attached to the relevant plan node.
Then $f$ is a partial elementary map of $M$.

Consequently, the complete type of a finite tuple over $C$ is determined by its
pure tree type over $C$ together with the finitely many local component types of
the active components.
\end{theorem}

\begin{proof}
Choose a block enumeration $C=A_0\subseteq A_1\subseteq\cdots\subseteq A_m=A$
as in Lemma~\ref{lem:block-enumeration}, and set $A_i'=f(A_i)$
for each $i\leq m$.  We prove by induction on $i$ that
$
f\restriction A_i:A_i\longrightarrow A_i'
$ 
is a partial elementary map.  The case $i=0$ is immediate, since $f$ fixes $C$
pointwise.

Assume the assertion for $i<m$.  Since $M$ is countable and
$\aleph_0$-categorical, finite partial elementary maps extend to automorphisms.
Let $g\in\Aut(M)$ extend $f\restriction A_i$.  Replacing $A'$ by $g^{-1}(A')$
and $f$ by $g^{-1}\circ f$, we may suppose that
$
A_i=A_i'
   \quad\text{and}\quad
   f\restriction A_i=\operatorname{id}_{A_i}.
$ 
This replacement preserves the hypotheses, because automorphisms preserve
both the pure tree structure and the canonical component structures.

Let $e=e_{i+1}$ and let $e'=f(e)$.  The pure tree isomorphism $f$ induces an
isomorphism
$$
\iota:\tcl(A_i e)\longrightarrow \tcl(A_i e')
$$
fixing $A_i$ pointwise and sending $e$ to $e'$.  Since
$\operatorname{pred}(e)\in A_i$ and $f$ fixes $A_i$ pointwise, we have
$\operatorname{pred}(e)=\operatorname{pred}(e').$
Also $\pi(e)=\pi(e')\in I(\Gamma)$.  The component hypothesis says precisely that
$e\equiv_{A_i}^{\loc} e'.$
By NI2, $\tp_M(\tcl(A_i e)/A_i)=\tp_M(\tcl(A_i e')/A_i),$
with the two blocks enumerated according to $\iota$.  Hence
$f\restriction A_{i+1}$ is partial elementary.  This completes the induction.

For the last assertion, apply what we have proved to
$ 
A=\tcl(C\bar a)
   \,\, \text{and}\,\,
   A'=\tcl(C\bar a')
$ 
for finite tuples $\bar a,\bar a'$. Their pure tree types and the local
types in the active components determine a map as in the first assertion.
This map is partial elementary and sends $\bar a$ to $\bar a'$ over $C$.
\end{proof}

\subsection{The tree-product characterization theorem}

Let $M$ be an admissible NIC structure. We write $P^{\mathrm{can}}(M)$
for the expansion of its underlying tree $\Gamma(S)$ by all canonical
component relations. More precisely, let $\sigma\in I(\Gamma)$,
let $R_p\in L_\sigma^{\mathrm{can}}$ be $n$-ary, and write
$\sigma=\tau\concat\langle k\rangle$. We interpret $\widehat R_p$ as follows:
$$
\begin{aligned}
\widehat R_p^{\,P^{\mathrm{can}}(M)}
=
\bigg\{&(b,b\concat\langle(k,s_0)\rangle,\ldots,
      b\concat\langle(k,s_{n-1})\rangle):\\
&\pi(b)=\tau
\text{ and }
(b,b\concat\langle(k,s_0)\rangle,\ldots,
   b\concat\langle(k,s_{n-1})\rangle)\models p\bigg\}.
\end{aligned}
$$
Thus $P^{\mathrm{can}}(M)$ is the canonical tree product obtained from the
canonical component languages of $M$.

\begin{theorem}
\label{thm:characterization}
Let $M$ be a countable admissible NIC expansion of $\Gamma(S)$, and let
$P^{\mathrm{can}}(M)$ be its canonical tree product.  Then $M$ and
$P^{\mathrm{can}}(M)$ have the same parameter-definable subsets of finite powers
of the underlying tree.  Equivalently, after passing to the canonical component
language, $M$ is interdefinable with its canonical tree product.

In particular, the complete theory of the canonical tree-product expansion is
determined by the tree plan together with the complete theories of the corresponding
canonical components.
\end{theorem}

\begin{proof}
First, every canonical component relation comes from a complete
$\emptyset$-type of $M$. These types are isolated by
$\aleph_0$-categoricity. Therefore the canonical component relations,
and hence all relations of $P^{\mathrm{can}}(M)$, are definable in $M$.

Conversely, let $D\subseteq M^r$ be definable in $M$ with parameters from a
finite set $E$.  Let 
$C=\tcl(E).$
The set $C$ is finite and tree-closed, and $E\subseteq C\subseteq\dcl(E)$
in the pure tree language. It is therefore enough to show that each
$C$-definable subset of $M^r$ is definable over $C$ in
$P^{\mathrm{can}}(M)$.

By $\aleph_0$-categoricity, there are only finitely many complete $r$-types
over $C$, and $D$ is a union of some of them. Theorem~\ref{thm:normal-form}
says that each type is determined by its pure tree type over $C$ and its
local component types. The pure tree type is definable in the
$L_\Gamma$-reduct. The local types are named in the canonical component
language. Thus each of these complete types, and consequently $D$, is
definable over $C$ in $P^{\mathrm{can}}(M)$.

The same argument proves the last assertion, since the finite tree plan
and the complete component theories determine all finite types in the
canonical expansion.
\end{proof}

\begin{remark}
By Theorem~\ref{thm:characterization}, admissible nil-interaction structures
are, up to definable expansion, canonical tree products. This explains
why the generic bipartite graph, presented as two sibling components,
is excluded: its edge relation introduces an additional relation between
the components.
\end{remark}

\subsection{Existence of admissible NIC structures}

\begin{proposition}\label{prop:tree-product-nic}
Let $\Gamma$ be a tree plan.  For each $\sigma\in I(\Gamma)$, let $T_\sigma$ be an
$\aleph_0$-categorical complete theory in a relational language $L_\sigma$.
Assume that $T_\sigma$ eliminates quantifiers, has trivial algebraic closure, and
has a unique complete $1$-type over $\emptyset$.  Let $N_\sigma$ be the countable
model of $T_\sigma$, and let
$
M=\prod_{\sigma\in I(\Gamma)}^{\!*} N_\sigma.
$
Then $M$ is $\aleph_0$-categorical, eliminates quantifiers in the tree-product
language, satisfies $\acl_M=\tcl$, and is an admissible NIC structure.
\end{proposition}

\begin{proof}
We first show that $M$ is $\aleph_0$-categorical. A finite tuple meets
finitely many components. Its quantifier-free type is described by its
pure tree diagram and the quantifier-free diagrams in those components.
Since the tree plan is finite and each component theory is
$\aleph_0$-categorical, there are only finitely many possibilities for
a fixed tuple length. Hence $M$ is $\aleph_0$-categorical.

We now show quantifier elimination. Let a finite partial isomorphism be
given. It preserves the pure tree diagram and the quantifier-free diagrams
in the active components. By quantifier elimination in the components,
its local restrictions are partial elementary. We extend the map by a
back-and-forth argument, adding one component element at a time and then
its forced $1$-successors. There are only finitely many pure tree diagrams
of a fixed arity, and the component extensions can be made independently.
Thus tuples with the same quantifier-free type have the same complete type.

We check algebraic closure.  Clearly $\tcl(A)\subseteq\dcl_M(A)\subseteq\acl_M(A),$
since predecessors, meets, and forced $1$-successors are definable. For
the reverse inclusion, let $A$ be finite and tree-closed, and let $b\notin A$.
Consider the first new non-forced point on the path from $[b\wedge A]$ to
$b$. It belongs to an infinite component. Triviality of algebraic closure
in that component gives infinitely many realizations of its local type
over the parameters from $A$. We move this point and copy the finite
subtree above it. Since the components are independent, each such local
move extends to an automorphism fixing $A$. This gives infinitely many
conjugates of $b$ over $A$, so $b\notin\acl_M(A)$.

NI1 follows from the unique component $1$-types and the symmetries of
the pure tree. For NI2, the relations on a new block are determined by
its tree position and the local type of its new infinity element.
This follows directly from the construction. Hence $M$ is admissible NIC.
\end{proof}

\subsection{Examples and boundary cases}\label{subsec:examples}

\begin{example}
Suppose that each $T_\sigma$ is the theory of an infinite pure set.
Then the product is the pure tree $\Gamma(S)$ with the predicates
$P_\tau$. In the measuring construction below, the dimension of a node
is the number of infinity-nodes on its path from the root.
\end{example}

\begin{example}\label{ex:random-components}
Let $T_\sigma$ be the theory of the countable random graph for every
$\sigma\in I(\Gamma)$.  The random graph is $\aleph_0$-categorical, eliminates
quantifiers in the graph language, has trivial algebraic closure, has a unique
complete $1$-type, and is $\MS$-measurable.  Hence the tree product of independent
random-graph components is an admissible NIC structure and is $\MS$-measurable
by Theorem~\ref{thm:main} below.
\end{example}

\begin{example}
\label{ex:vector-space-boundary}
Let $V$ be an infinite-dimensional vector space over a fixed finite field.
Its theory is $\aleph_0$-categorical and measurable. However, algebraic
closure is linear span, so the component closure is not trivial. To include
tree products with these components, one would need to replace $\tcl$ by
a closure which also takes local linear spans. The finite-dimensional
approximations of these spaces fit the multidimensional asymptotic setting
of Section~\ref{sec:mac-trees}.
\end{example}

\section{A type criterion for \texorpdfstring{$\MS$}{MS}-measurability}
\label{sec:omega-cat-criterion}

We recall $\MS$-measurability and give a criterion using complete types over
finite parameter sets. We assign the empty set the value $(0,0)$ and take
$0\in\N$. By the Ryll--Nardzewski theorem, complete types over finite
parameter sets in an $\aleph_0$-categorical structure are isolated; see
\cite{Marker2002}.

\begin{definition}
Let $M$ be an infinite structure.  Write $\Def^n(M)$ for the collection of
parameter-definable subsets of $M^n$, and let
$ 
\Def(M)=\bigcup_{n\geq1}\Def^n(M).
$ 
The structure $M$ is \emph{$\MS$-measurable} if there is a function
$$
h=(\dimn,\meas):\Def(M)\longrightarrow
   (\N\times\R_{>0})\cup\{(0,0)\}
$$
satisfying the following axioms.

\begin{enumerate}[label=\textup{MS\arabic*.},leftmargin=4.0em]
\item For every formula $\varphi(\bar x,\bar y)$, the set $\{h(\varphi(M^{|\bar x|},\bar a)):\bar a\in M^{|\bar y|}\}$
is finite.

\item If $X$ is finite, then
$h(X)=(0,|X|).$

\item For every formula $\varphi(\bar x,\bar y)$ and every pair $(d,\mu)$
occurring among the values
$$
h(\varphi(M^{|\bar x|},\bar a)),
   \qquad \bar a\in M^{|\bar y|},
$$
the parameter set
$$
\{\bar a\in M^{|\bar y|}:
      h(\varphi(M^{|\bar x|},\bar a))=(d,\mu)\}
$$
is $\emptyset$-definable.

\item \textup{(Fubini)}
Let $f:X\to Y$ be a definable surjection.  For each pair $(d,\mu)$ occurring
as the value of a fiber, let
$Y_{d,\mu}=\{y\in Y:h(f^{-1}(y))=(d,\mu)\}.$
Then the nonempty $Y_{d,\mu}$ form a finite definable partition of $Y$, and $\dimn(X)=\max_{(d,\mu)}\{d+\dimn(Y_{d,\mu})\},$
while
$$
\meas(X)=
   \sum_{d+\dimn(Y_{d,\mu})=\dimn(X)}
   \mu\cdot\meas(Y_{d,\mu}).
$$
\end{enumerate}
\end{definition}

\begin{lemma}\label{lem:finite-unions}
Let $M$ be $\MS$-measurable, and let $X_0,\ldots,X_{k-1}$
be pairwise disjoint definable subsets of the same Cartesian power of $M$.  Then
$$
\dimn\Big(\bigsqcup_{i<k}X_i\Big)=\max_{i<k}\dimn(X_i),
$$
and
$$
\meas\Big(\bigsqcup_{i<k}X_i\Big)=
   \sum_{\dimn(X_i)=\max_j\dimn(X_j)}\meas(X_i).
$$
\end{lemma}

\begin{proof}
If every $X_i$ is empty, the result follows from $h(\emptyset)=(0,0)$.
Otherwise, we discard the empty sets and assume that each $X_i$ is
nonempty. Let
$
X=\bigsqcup_{i<k}X_i.
$ 
Choose distinct elements $c_0,\ldots,c_{k-1}\in M$ and let 
$Y=\{c_0,\ldots,c_{k-1}\}.$
Define a definable surjection $f:X\to Y$ by
$$
f(x)=c_i \quad\text{if and only if}\quad x\in X_i.
$$
By MS2, each singleton $\{c_i\}$ has dimension $0$ and measure $1$.  Applying
MS4 gives $\dimn(X)=\max_{i<k}\dimn(X_i)$
and
$$
\meas(X)=
   \sum_{\dimn(X_i)=\max_j\dimn(X_j)}\meas(X_i).
$$
\end{proof}

For an $\aleph_0$-categorical structure $M$, finite $C\subseteq M$, and
$a\in M^n$, write $h(a/C)$ for the value assigned to the definable set
$\tp(a/C)(M)$ of realizations of $\tp(a/C)$. Lemma~\ref{lem:finite-unions}
gives the value of a finite disjoint union of such types.

\begin{theorem}\label{thm:omega-cat-criterion}
Let $M$ be a countably infinite $\aleph_0$-categorical structure.  Then $M$ is
$\MS$-measurable if and only if there is an assignment $(a,C)\longmapsto h(a/C)=(\dimn(a/C),\meas(a/C))$
for all finite $C\subseteq M$ and all finite tuples $a$ from $M$, satisfying the
following conditions.

\begin{enumerate}[label=\textup{cMS\arabic*.},leftmargin=4.2em]
\item \textup{(Automorphism invariance)}
For every $g\in\Aut(M)$,
$h(ga/gC)=h(a/C).$

\item \textup{(Finite unions of complete types)}
If $X\subseteq M^r$ is $C$-definable and $X=p_0(M)\dsum\cdots\dsum p_{k-1}(M),$
where each $p_i\in S_r(C)$ is a complete type over $C$, then, writing
$$
h(p_i)=h(a_i/C)
   \quad\text{for any }a_i\models p_i,
$$
we have $\dimn(X)=\max_{i<k}\dimn(p_i),$
and
$$
\meas(X)=
   \sum_{\dimn(p_i)=\dimn(X)}\meas(p_i).
$$
The empty definable set is assigned value $(0,0)$.  This rule is required to be
independent of the finite parameter set over which $X$ is decomposed.

\item \textup{(Algebraic types)}
If $a\in\acl(C)^n$, then
$$
\dimn(a/C)=0,
   \qquad
   \meas(a/C)=|\{ga:g\in\Aut(M/C)\}|.
$$

\item \textup{(Product over definable dependence)}
If $b\in\dcl(Ca)^m$, then
$$
\begin{aligned}
   \dimn(a/C)&=\dimn(b/C)+\dimn(a/Cb),\\
   \meas(a/C)&=\meas(b/C)\cdot\meas(a/Cb).
\end{aligned}
$$
\end{enumerate}
\end{theorem}

\begin{proof}
Suppose first that $h$ satisfies cMS1--cMS4. Let $X\subseteq M^r$ be a
nonempty set definable over a finite set $C$. We write $X$ as a disjoint
union of complete $r$-types over $C$ and define $h(X)$ by cMS2.
We set $h(\emptyset)=(0,0)$. The independence requirement in cMS2 shows
that this definition does not depend on the choice of $C$.

We verify the Macpherson--Steinhorn axioms.

For MS1, fix a formula $\varphi(\bar x,\bar y)$. There are only finitely
many complete types of $\bar y$ over $\emptyset$. Suppose that $\bar a$
and $\bar a'$ have the same type. Choose an automorphism sending
$\bar a$ to $\bar a'$. By cMS1, corresponding complete types in the
decompositions of
$
\varphi(M^{|\bar x|},\bar a)
   \,\, \text{and}\,\,
   \varphi(M^{|\bar x|},\bar a')
$ 
have the same values. Applying cMS2, the two definable sets have the
same value as well. Thus only finitely many values occur.

MS2 follows from cMS2 and cMS3, since a finite definable set is a finite disjoint
union of algebraic complete types.

For MS3, choose a formula $\theta_q(\bar y)$ isolating each complete
type $q(\bar y)$ over $\emptyset$. By cMS1, the set of parameters
$\bar a$ satisfying $h(\varphi(M^{|\bar x|},\bar a))=(d,\mu)$
is the union of the finitely many $\theta_q(M)$ with corresponding value
$(d,\mu)$. Hence it is $\emptyset$-definable.

It remains to prove MS4.  Let $f:X\to Y$
be a $C$-definable surjection.  If $Y=\emptyset$, then $X=\emptyset$ and the
assertion is immediate.  Otherwise, decompose $Y=q_0(M)\dsum\cdots\dsum q_{r-1}(M)$
into complete types over $C$, and choose $b_i\models q_i$.  For fixed $i$,
decompose the fiber $f^{-1}(b_i)=p_{i,0}(M)\dsum\cdots\dsum p_{i,k_i-1}(M)$
into complete types over $Cb_i$.  Choose $a_{i,j}\models p_{i,j}$ and let
$X_{i,j} = \{a\in X:\tp(a,f(a)/C)=\tp(a_{i,j},b_i/C)\}.$
Each $X_{i,j}$ is the realization set of the complete type
$\tp(a_{i,j}/C).$
Indeed, $f$ is $C$-definable, so $\tp(a/C)$ determines
$\tp(a,f(a)/C)$, and the converse is immediate. Therefore the sets
$X_{i,j}$ form a finite disjoint partition of $X$ into complete types
over $C$.

Since $f(a_{i,j})=b_i$, we have
$b_i\in\dcl(Ca_{i,j}).$
Therefore cMS4 gives
$$
\begin{aligned}
   \dimn(a_{i,j}/C)&=\dimn(b_i/C)+\dimn(a_{i,j}/Cb_i),\\
   \meas(a_{i,j}/C)&=\meas(b_i/C)\cdot\meas(a_{i,j}/Cb_i).
\end{aligned}
$$

For each $i$, we apply cMS2 to the decomposition of $f^{-1}(b_i)$ into
the types $p_{i,j}$. By cMS1, the fibers over realizations of $q_i$ all
have the same value. Write
$$
h(f^{-1}(b_i))=(d_i,\mu_i),
   \qquad
   h(b_i/C)=(e_i,\nu_i).
$$
Then
$$
Y_{d,\mu}
   =
   \bigsqcup_{\{i:(d_i,\mu_i)=(d,\mu)\}} q_i(M),
$$
so the sets $Y_{d,\mu}$ are finite unions of complete types over $C$ and are
therefore definable.  Applying cMS2 to the partition of $X$ by the sets
$X_{i,j}$ gives
$\dimn(X)=\max_i(d_i+e_i).$
Only the terms with $d_i+e_i=\dimn(X)$ contribute to the measure.
Their sum is
$$
\sum_{d_i+e_i=\dimn(X)} \mu_i\nu_i.
$$ 
This is precisely the Fubini formula
$$
\meas(X)=
   \sum_{d+\dimn(Y_{d,\mu})=\dimn(X)}
   \mu\cdot\meas(Y_{d,\mu}).
$$
Thus MS4 holds.

Conversely, suppose $M$ is $\MS$-measurable.  Define $h(a/C)$
to be the value of the definable set $\tp(a/C)(M)$ of realizations of
$\tp(a/C)$.

We verify cMS1--cMS4.  For cMS1, let $g\in\Aut(M)$.  If
$\theta(x,c)$ isolates $\tp(a/C)$, where $c$ enumerates $C$, then
$\theta(x,gc)$ isolates $\tp(ga/gC)$.  By MS3, the parameters $c$ and
$gc$ lie in the same $\emptyset$-definable part of the relevant
dimension and measure partition.  Hence
$h(ga/gC)=h(a/C).$

cMS2 follows from Lemma~\ref{lem:finite-unions}. The value is independent
of the parameter set because the original measuring function is defined
on definable sets.

cMS3 follows from MS2 and $\aleph_0$-categoricity.  If $a\in\acl(C)^n$, then the
realizations of $\tp(a/C)$ are exactly the finite orbit
$\{ga:g\in\Aut(M/C)\}.$
Therefore MS2 gives
$$
\dimn(a/C)=0,
   \qquad
   \meas(a/C)=|\{ga:g\in\Aut(M/C)\}|.
$$

It remains to prove cMS4.  Suppose $b\in\dcl(Ca)^m$.  Since $\tp(a/C)$ is
isolated, there is a $C$-definable function $F$, defined on the realization set $X=\tp(a/C)(M),$
such that
$F(a)=b.$
Let
$Y=\tp(b/C)(M).$
By homogeneity over the finite set $C$, the map $F:X\to Y$
is surjective.

We claim that the fiber over $b$ is the realization set of $\tp(a/Cb)$.
On this fiber, a formula $\psi(x,b,C)$ is equivalent to the $C$-formula
$\psi(x,F(x),C)$. Therefore
$F^{-1}(b)=\tp(a/Cb)(M).$
By automorphism invariance over $C$ and MS3, all fibers over $Y$ have
the same value. Applying MS4 to $F:X\to Y$
gives $\dimn(a/C)=\dimn(b/C)+\dimn(a/Cb)$
and
$\meas(a/C)=\meas(b/C)\cdot\meas(a/Cb).$
Thus cMS4 holds.
\end{proof}

\section{Lifting \texorpdfstring{$\MS$}{MS}-measurability}
\label{sec:lifting}

Throughout this section, let $M$ be an admissible NIC expansion of
$\Gamma(S)$. For each $\sigma\in I(\Gamma)$, we write $T_\sigma$ for the
common theory of the canonical components $M[e]$ with $\pi(e)=\sigma$.
By the normal form, components of the same plan-type are elementarily
equivalent in their common canonical language.

We first note that each $T_\sigma$ is $\aleph_0$-categorical. Indeed,
component $n$-types are coded by complete types in $M$ of tuples $(b,a_0,\ldots,a_{n-1}),$
where $b$ is the predecessor and $a_0,\ldots,a_{n-1}$ belong to the
component of plan-type $\sigma$ above $b$. There are only finitely many
such types, by $\aleph_0$-categoricity of $M$. In particular, complete
types over finite parameter sets in a component are isolated.

We also need that algebraic closure in each canonical component is
trivial. Let $D\subseteq\Comp(e)$ be finite, and suppose that
$d\in\acl_{M[e]}(D)$. Let $b=\operatorname{pred}(e)$. By admissibility,
$d\in\acl_M(Db)=\tcl(Db).$
Intersecting $\tcl(Db)$ with $\Comp(e)$ gives $D$.  Hence $d\in D$.

Assume that each $T_\sigma$ is $\MS$-measurable.  Fix, for each
$\sigma\in I(\Gamma)$, a measuring function $h_\sigma=(\dimn_\sigma,\meas_\sigma)$
with
$\dimn_\sigma(x=x)=1.$
We use this normalization to simplify the notation. Without it, the same
proof gives the corresponding weighted sum of local component dimensions.

For the empty tuple in a component we use the convention
$h_\sigma(\emptyset/D)=(0,1).$

\subsection{The global type assignment}

Let $C\subseteq M$ be finite and tree-closed, and let $a$ be a finite tuple of
elements of $M$.  Let
$A=\tcl(Ca).$
For every active component $K$ in $A/C$, let
$ 
A_K=A\cap K,
   \,\, 
   C_K=C\cap K,
$ 
and let $\bar a_K$ enumerate $A_K\setminus C_K$.  If $K$ has plan node
$\sigma$, define $h_K(A/C)=h_\sigma(\bar a_K/C_K),$
where the value on the right is computed in the canonical component of
plan-type $\sigma$. It does not depend on the enumeration, since
permutations of coordinates are definable bijections.

\begin{definition}
For a complete type $\tp(a/C)$ with $C$ tree-closed, define
$$
\dimn(a/C)=\sum_{K\in\mathcal K(A/C)}\dimn_K(A/C),
$$
and
$$
\meas(a/C)=\prod_{K\in\mathcal K(A/C)}\meas_K(A/C),
$$
where $\mathcal K(A/C)$ is the finite set of active components in $A/C$.  The
empty sum is $0$, and the empty product is $1$.

For a finite set $C$ which is not necessarily tree-closed, set
$h(a/C)=h(a/\tcl(C)).$
For an arbitrary definable set, we extend $h$ by finite disjoint unions of
complete types as in cMS2 of Theorem~\ref{thm:omega-cat-criterion}.
\end{definition}

If $C$ is tree-closed and $[e\wedge C]=e_0<e_1<\cdots<e_m=e$
is the path from $C$ to $e$, then
$$
\dimn(e/C)=|\{i:1\leq i\leq m\text{ and }\pi(e_i)\in I(\Gamma)\}|.
$$
To see this, note that the new infinity-nodes are precisely the points on
the path which contribute active components. By normalization and local
algebraic triviality, each contributes dimension $1$. The forced
$1$-successors contribute dimension $0$.

If $\pi(e)=\sigma\in I(\Gamma)$ and $\operatorname{pred}(e)\in C$, then
$h(e/C)=h_\sigma(e/C\cap\Comp(e)).$

\subsection{Base change}

\begin{lemma}\label{lem:base-change}
Let $B\subseteq C\subseteq M$ be finite tree-closed sets, and let
$p(\bar x)\in S_{\bar x}(B)$ be a complete type.  Write the set of realizations of
$p$ as a finite disjoint union of complete types over $C$,
$p(M)=p_0(M)\dsum\cdots\dsum p_{r-1}(M).$
Then $\dimn(p)=\max_{i<r}\dimn(p_i),$
and
$$
\meas(p)=
   \sum_{\dimn(p_i)=\dimn(p)}\meas(p_i).
$$
\end{lemma}

\begin{proof}
We first note that, in a component, every nonalgebraic $1$-type over a
finite set has dimension $1$. Indeed, the universe has dimension $1$,
and a definable set of dimension $0$ is finite. For the latter assertion,
an infinite set of dimension $0$ would contain arbitrarily many disjoint
singletons, contrary to the finite-union rule and finiteness of its
measure. Since algebraic closure in each component is trivial, successive
coordinate projections and Fubini show that a tuple of $s$ distinct
elements outside a finite parameter set has dimension $s$. Hence, for
every finite tree-closed set $D$ and finite tuple $\bar a$,
\[
\dimn(\bar a/D)
=\big|\{u\in\tcl(D\bar a)\setminus D:
                  \pi(u)\in I(\Gamma)\}\big|.
\tag{*}
\]

We prove the assertion first when $C=\tcl(Bc)$, where
$c\notin B$, $\pi(c)\in I(\Gamma)$, and
$\operatorname{pred}(c)\in B$. Let $K=\Comp(c)$,
$B_K=B\cap K$, and $d=\dimn(p)$. By
Theorem~\ref{thm:characterization}, we may work in the canonical tree
product.

Let $\bar a\models p$ and $A=\tcl(B\bar a)$. Tree closure commutes
with unions, so $\tcl(C\bar a)=C\cup A$. Also, $c$ is the only element
of $C\setminus B$ whose plan node belongs to $I(\Gamma)$. By $(*)$,
\[
\dimn(\bar a/C)=
\begin{cases}
d,&c\notin A,\\
d-1,&c\in A.
\end{cases}
\]
Thus it remains to find the extensions for which $c\notin A$ and to
compute the sum of their measures.

Choose a realization $\bar a\models p$ and enumerate
$(A\cap K)\setminus B_K$ by a tuple $\bar u$ of distinct elements.
The entries of $\bar u$ are ancestors of entries of $\bar a$. We use
the same ancestor terms to obtain the corresponding tuple from every
other realization of $p$. Let
$
q=\tp_K(\bar u/B_K),\text{ with } s=|\bar u|.
$
Here $\tp_K$ and $h_K$ are computed in the canonical component
structure on $K$. If $s=0$, we use the empty tuple and its value
$(0,1)$.

Let $q_0,\ldots,q_{t-1}$ be the distinct complete extensions of $q$
over $B_Kc$, and let
\[
J=\{j<t:q_j\text{ implies }u_\ell\ne c
                         \text{ for every }\ell<s\}.
\]
By the first paragraph, $\dimn_K(q)=s$, and
$\dimn_K(q_j)=s$ exactly when $j\in J$. An extension identifying an
entry with $c$ has dimension $s-1$. If $s=0$, there is just one
extension, with value $(0,1)$. For $s>0$, apply the finite-union rule
in $K$. In either case we obtain
\[
J\ne\emptyset,\qquad
\meas_K(q)=\sum_{j\in J}\meas_K(q_j).
\tag{**}
\]

We claim that the extensions of $p$ over $C$ of dimension $d$
correspond bijectively to the types $q_j$, $j\in J$. To see existence,
choose $\bar u'\models q_j$. Since $K$ is countable and
$\aleph_0$-categorical, there is an automorphism of $K$ over $B_K$
sending $\bar u$ to $\bar u'$. This automorphism extends to an
automorphism of the tree product over $B$: move the subtrees above the
points of $K$ accordingly, and leave the subtrees above $B_K$ fixed.
This fixes $B$, since $B$ is tree-closed. The image $\bar a'$ of
$\bar a$ realizes $p$, and $c\notin\tcl(B\bar a')$ because no entry
of $\bar u'$ is $c$.

For uniqueness, suppose that two realizations of $p$ give the same
$q_j$, with $j\in J$. For each such realization $\bar a'$,
$\tcl(B\bar a')\cap C=B$.
Indeed, a point in $(C\setminus B)\cap\tcl(B\bar a')$ would have
$c$ as an ancestor. The only additional local parameter is therefore
$c$ in $K$; all other local types are those already specified by $p$.
Theorem~\ref{thm:normal-form} now shows that the two realizations have
the same type over $C$. Conversely, equality of their types over $C$
implies equality of the corresponding local types over $B_Kc$.
This proves the claim.

Let $\nu$ be the product of the local measures defining $\meas(p)$,
with the factor from $K$ omitted. If $K$ is not active, take its
omitted factor to be $1$. By construction,
$
\meas(p)=\nu\meas_K(q),
$ 
and the extension corresponding to $q_j$ has measure
$\nu\meas_K(q_j)$. Therefore $(**)$ gives
\[
\meas(p)
=\nu\sum_{j\in J}\meas_K(q_j)
=\sum_{\dimn(p_i)=d}\meas(p_i).
\]
Since $J\ne\emptyset$, the maximum dimension of an extension is $d$.
This proves the assertion for one block.

For general $B\subseteq C$, choose a block enumeration from $B$ to $C$
as in Lemma~\ref{lem:block-enumeration}, and apply the preceding
argument at each step. Dimensions never increase. Thus a final
extension has dimension $\dimn(p)$ exactly when its dimension is
preserved at every step. Summing the corresponding measure identities
successively counts each such final extension exactly once and gives
the required equality.
\end{proof}

\begin{corollary}\label{cor:well-defined}
The value assigned to a definable set by decomposing it into complete types is
independent of the finite parameter set used to define it.
\end{corollary}

\begin{proof}
Let $X$ be definable over finite parameter sets $E_1$ and $E_2$.  Let 
$
B_i=\tcl(E_i)
   \,\,\text{for }i=1,2,
$
and let
$C=\tcl(B_1B_2).$
For each $i$, refine the types over $B_i$ occurring in $X$ to types over
$C$. By Lemma~\ref{lem:base-change}, this does not change the computed
value. Both decompositions therefore give the value obtained over $C$.
\end{proof}

\subsection{Verification of the criterion}

\begin{lemma}\label{lem:cMS1-global}
The global assignment satisfies cMS1.
\end{lemma}

\begin{proof}
Let $g\in\Aut(M)$. Since $g$ preserves tree closure and the canonical
component structures, it sends the active components for $a/C$ to those
for $ga/gC$, preserving their local types. Each component measuring
function is automorphism invariant. Hence
$h(ga/gC)=h(a/C).$
\end{proof}

\begin{lemma}\label{lem:cMS2-global}
The global assignment satisfies cMS2.
\end{lemma}

\begin{proof}
The finite-union rule is part of the definition of $h$.
Its independence of the parameter set is
Corollary~\ref{cor:well-defined}.
\end{proof}

\begin{lemma}\label{lem:cMS3-global}
The global assignment satisfies cMS3.
\end{lemma}

\begin{proof}
If $a\in\acl(C)^n$, then $a\in\tcl(C)^n$ by admissibility.  Since the functions
defining tree closure are in the language, we have
$\tcl(C)\subseteq\dcl(C).$
Hence $a$ has a singleton orbit over $C$. Also, $a$ adds no active
component over $\tcl(C)$. By the conventions for the empty sum and
product,
$
\dimn(a/C)=0
   \,\,\text{and }\, 
   \meas(a/C)=1,
$ 
as required by cMS3.
\end{proof}

\begin{lemma}\label{lem:cMS4-global}
The global assignment satisfies cMS4.
\end{lemma}

\begin{proof}
Let $C$ be finite, let $a$ and $b$ be finite tuples, and assume
$b\in\dcl(Ca)^m.$
Replace $C$ by $C_0=\tcl(C)$
and let 
$
B=\tcl(C_0b),
   \,\,
   A=\tcl(C_0a).
$
Since $b\in\dcl(Ca)\subseteq\acl(Ca)=\tcl(Ca)=\tcl(C_0a)$
by admissibility, we have
$A=\tcl(C_0ab).$
Thus
$C_0\subseteq B\subseteq A.$

For each component $K$ active in $A/C_0$, write
$
C_{0,K}=C_0\cap K,
   \,\, 
   B_K=B\cap K,
   \,\, 
   A_K=A\cap K.
$
Let $\bar a_K$ enumerate $A_K\setminus C_{0,K}$ and let $\bar b_K$ enumerate
$B_K\setminus C_{0,K}$.  Up to a permutation, $\bar b_K$ is a subtuple of $\bar a_K$. Therefore $\bar b_K\in\dcl(C_{0,K}\bar a_K)$
in the canonical component structure. Applying cMS4 in that component,
we get $\dimn_K(A_K/C_{0,K}) = \dimn_K(B_K/C_{0,K})+\dimn_K(A_K/B_K)$
and
$$
\meas_K(A_K/C_{0,K})
   =
   \meas_K(B_K/C_{0,K})\cdot\meas_K(A_K/B_K).
$$
If one of the relevant local tuples is empty, we use the convention
$h_\sigma(\emptyset/D)=(0,1).$

We now sum the dimension equalities and multiply the measure equalities
over the active components in $A/C_0$. This gives $\dimn(a/C)=\dimn(b/C)+\dimn(a/Cb)$
and
$\meas(a/C)=\meas(b/C)\cdot\meas(a/Cb).$
Thus cMS4 holds.
\end{proof}

\begin{theorem}[Lifting $\MS$-measurability]\label{thm:main}
Let $\Gamma$ be a tree plan and let $M$ be an admissible NIC expansion of
$\Gamma(S)$.  Suppose that for every $\sigma\in I(\Gamma)$ the component theory
$T_\sigma$ is $\MS$-measurable via a measuring function $h_\sigma=(\dimn_\sigma,\meas_\sigma)$
with
$\dimn_\sigma(x=x)=1.$
Then $M$ is $\MS$-measurable.

More precisely, let $C$ be finite and tree-closed, and let $a$ be a
tuple. To compute the value of $\tp(a/C)$, take the active components
of $\tcl(Ca)$ over $C$, add their local dimensions, and multiply their
local measures.
\end{theorem}

\begin{proof}
By Lemmas~\ref{lem:cMS1-global}, \ref{lem:cMS2-global},
\ref{lem:cMS3-global}, and \ref{lem:cMS4-global}, the function constructed
above satisfies cMS1--cMS4. Theorem~\ref{thm:omega-cat-criterion} now
shows that $M$ is $\MS$-measurable.
\end{proof}

\begin{corollary}\label{cor:tree-product-measurable}
Let $\Gamma$ be a tree plan.  For each $\sigma\in I(\Gamma)$, let $T_\sigma$ be
$\aleph_0$-categorical, algebraically trivial, quantifier-eliminating, and have a
unique complete $1$-type over $\emptyset$.  Suppose moreover that $T_\sigma$ is
$\MS$-measurable and normalized with
$\dimn_\sigma(x=x)=1.$
If $N_\sigma$ is the countable model of $T_\sigma$, then the canonical
$\Gamma$-tree product
$ 
\prod_{\sigma\in I(\Gamma)}^{\!*}N_\sigma
$ 
is $\MS$-measurable.
\end{corollary}

\begin{proof}
By Proposition~\ref{prop:tree-product-nic}, the tree product is admissible NIC.
Apply Theorem~\ref{thm:main}.
\end{proof}

\begin{corollary}\label{cor:dimension-tree}
For $\tau\in\Gamma$, let
$w(\tau)=|\{\sigma\in I(\Gamma):\sigma\leq\tau\}|.$
Then in the measurable structure of Theorem~\ref{thm:main}, $\dimn(P_\tau)=w(\tau),$
and
$\dimn(M)=\max_{\tau\in\Gamma}w(\tau).$
\end{corollary}

\begin{proof}
To obtain a generic element of $P_\tau$, we make an independent choice
in each infinity component on the path from the root to $\tau$.
Each choice contributes dimension $1$, and each forced $1$-successor
contributes dimension $0$. Therefore
$\dimn(P_\tau)=w(\tau).$
The statement for $M$ follows because
$
M=\bigcup_{\tau\in\Gamma}P_\tau
$
is a finite union.
\end{proof}

\section{Rank, one-basedness, and tree envelopes}
\label{sec:rank-envelopes}

We prove preservation of finite $\SU$-rank, one-basedness, and smooth
approximability for canonical tree products.

\subsection{A coordinate presentation}

For $\tau\in\Gamma$, let
$I_\tau=\{\sigma\in I(\Gamma):\sigma\leq\tau\}.$
Fix once and for all an ordering of each finite set $I_\tau$ compatible with the
tree order.

Let
$
M=\prod_{\sigma\in I(\Gamma)}^{\!*}N_\sigma
$ 
be the canonical $\Gamma$-tree product.  If $a\in P_\tau(M)$ and
$\sigma\in I_\tau$, let $c_\sigma(a)$ denote the coordinate of the unique ancestor
of $a$ having plan node $\sigma$.  Thus, if that ancestor has the form $b\concat\langle(k,s)\rangle,$
then
$c_\sigma(a)=s.$

\begin{lemma}\label{lem:coordinate-presentation}
Let $D_\Gamma$ be the finite many-sorted disjoint union of the component
structures $N_\sigma$, with one sort for each $\sigma\in I(\Gamma)$.  After adding
the finitely many product sorts
$
\prod_{\sigma\in I_\tau}N_\sigma
   \,\, (\text{ for }\tau\in\Gamma),
$
where the empty product is a singleton, the canonical tree product $M$ is
definable without quotienting in this finite product-sort expansion of
$D_\Gamma$.

More precisely, for each $\tau\in\Gamma$, the predicate $P_\tau(M)$ is identified
with
$
\prod_{\sigma\in I_\tau}N_\sigma.
$
Under these identifications, the tree order, meet, predecessor, and the predicates
$P_\tau$ are definable from equality, coordinate projections, and the finite plan
$\Gamma$.  If $R\in L_\sigma$ and $\sigma=\tau\concat\langle k\rangle$, then the
corresponding relation $\widehat R$ is the pullback of $R$ along the
$\sigma$-coordinate, restricted to tuples whose coordinates below $\tau$ agree.
\end{lemma}

\begin{proof}
For $a\in P_\tau(M)$, define
$\rho_\tau(a)=(c_\sigma(a))_{\sigma\in I_\tau}.$
The map $\rho_\tau$ is a bijection onto the displayed product. If the
last step from $\operatorname{pred}(\tau)$ to $\tau$ is labelled by
$\infty$, the predecessor map removes the final coordinate. If it is
labelled by $1$, the infinity-coordinates remain unchanged. Since the
plan is finite, we can define the order and meet by considering finitely
many cases and comparing coordinates on common initial segments.
Thus the pure tree structure is definable using equality and projections.

Now suppose $R\in L_\sigma$ is $n$-ary and
$\sigma=\tau\concat\langle k\rangle.$
By definition of the canonical tree product,
$$
\widehat R(b,b\concat\langle(k,s_0)\rangle,\ldots,
                b\concat\langle(k,s_{n-1})\rangle)
$$
holds if and only if
$N_\sigma\models R(s_0,\ldots,s_{n-1}).$
Thus the entries $b\concat\langle(k,s_i)\rangle$
agree with $b$ in the coordinates below $\tau$, and their
$\sigma$-coordinates satisfy $R$ in $N_\sigma$. This is definable
in the indicated expansion by product sorts.
\end{proof}

\subsection{Rank and one-basedness}

If $K$ is a component, write $\sigma(K)$ for its plan node; that is, $\sigma(K)=\pi(e)$
for any $e\in K$.

\begin{theorem}\label{thm:rank-one-based}
Let
$
M=\prod_{\sigma\in I(\Gamma)}^{\!*}N_\sigma
$
be a canonical $\Gamma$-tree product. Assume that each component theory
$T_\sigma=\Th(N_\sigma)$ is supersimple of finite $\SU$-rank. Then
$\Th(M)$ is supersimple of finite $\SU$-rank. 
Moreover, let $C\subseteq M$ be finite and tree-closed, let $a$ be a
finite tuple, and let $A=\tcl(Ca)$. For each component $K$ active in
$A/C$, let
$
A_K=A\cap K,\,\, C_K=C\cap K,
$
and let $\bar a_K$ enumerate $A_K\setminus C_K$. Then
\[
\SU_M(a/C)
=\sum_{K\in\mathcal K(A/C)}
       \SU_{T_{\sigma(K)}}(\bar a_K/C_K).
\]
In particular, for every $\tau\in\Gamma$,
\[
\SU_M(P_\tau)
=\sum_{\substack{\sigma\in I(\Gamma)\\\sigma\leq\tau}}
       \SU(T_\sigma),
\]
and
\[
\SU(M)
=\max_{\tau\in\Gamma}
  \sum_{\substack{\sigma\in I(\Gamma)\\\sigma\leq\tau}}
       \SU(T_\sigma).
\]
If each $T_\sigma$ is one-based, then $\Th(M)$ is one-based.

The same conclusions hold for any admissible NIC structure whose
canonical component theories satisfy the corresponding hypotheses.
\end{theorem}

\begin{proof}
By Lemma~\ref{lem:coordinate-presentation}, $M$ is interpretable in
the finite many-sorted disjoint union of the component structures,
after adding finitely many definable product sorts. The theory of
this disjoint union is simple, and simplicity is preserved under
interpretations. Thus $\Th(M)$ is simple. We work in a sufficiently
saturated elementary extension and continue to write $M$ for this
extension. Its components are saturated models of the corresponding
theories $T_\sigma$.

We first describe types and formulas with variables in one component.
For this purpose, name every parameter-free definable relation in each
component, and name the forced $1$-successor functions in the tree.
These are definitional expansions. In these languages a map between
finite tree-closed sets which preserves the tree and is elementary in
each component extends by one block: choose the new infinity element
with the required local type, and then add its forced $1$-successors.
Saturation of the components permits this choice in either direction.
The resulting back-and-forth shows that such maps are elementary.
In particular, this gives quantifier elimination in the expanded
tree-product language. No assumption on local algebraic closure is
needed for this argument.

Fix a component $K$ of plan-type $\sigma$, and let $b$ be its common
predecessor. If $D$ is tree-closed and contains $b$, then the type in
$M$ of a tuple from $K$ over $D$ is determined by its type in $K$ over
$D\cap K$. Indeed, the points of $K$ occurring as ancestors of
parameters in $D$ already belong to $D\cap K$. The remaining
components introduce no further relations on $K$.

We will also use the following consequence for formulas. Let $B$ be
tree-closed with $b\in B$, and let $\varphi(x,\bar d)$ be a formula,
where $x$ ranges over $K$. On $K$, this formula is equivalent to a
local formula $\psi(x,\bar c)$, with parameters from $B\cap K$ and
a finite tuple
$
\bar c\subseteq\tcl(B\bar d)\cap K.
$
The entries of $\bar c$ are obtained from $\bar d$ by fixed
$B$-definable ancestor maps. Moreover, the same formula and maps work
for every $\bar d'\equiv_B\bar d$. To see this, apply the quantifier
elimination just established. A tree relation involving $x$ and a
parameter reduces to equality with an ancestor in $K$, or has a truth
value independent of $x$. A local relation involving $x$ is either
a relation in $K$, or forces $x$ to equal such an ancestor. The
parameter-only conditions are fixed by $\tp(\bar d/B)$.

\smallskip
\noindent\emph{Claim.} If $B\subseteq D$ are small tree-closed sets,
$b\in B$, and $e\in K$, then
\[
\tp_M(e/D)\text{ forks over }B
\quad\Longleftrightarrow\quad
\tp_K(e/D\cap K)\text{ forks over }B\cap K.
\tag{*}
\]

Since both theories are simple, it is enough to compare dividing.
Suppose first that a local formula $\psi(x,\bar c)$ divides over
$B\cap K$. Let $(\bar c_j)_{j<\omega}$ be a witnessing
$(B\cap K)$-indiscernible sequence in $K$, with $\bar c_0=\bar c$.
By the description of types above, this sequence is also
$B$-indiscernible in $M$. Relativize $\psi$ to the $b$-definable
component $K$. The resulting formula divides over $B$ in $M$.

Conversely, suppose that $\varphi(x,\bar d)\in\tp_M(e/D)$ divides
over $B$, and let $(\bar d_j)_{j<\omega}$ be a witnessing
$B$-indiscernible sequence, with $\bar d_0=\bar d$. We may conjoin
$x\in K$ to the formula. Use the same local formula $\psi$ and the
same ancestor maps for each $\bar d_j$, and write $\bar c_j$ for the
resulting tuple in $K$. Then $(\bar c_j)_{j<\omega}$ is
$(B\cap K)$-indiscernible in the component. Also,
$
\{\psi(x,\bar c_j):j<\omega\}
$ 
is inconsistent, since on $K$ these formulas are equivalent to the
corresponding instances of $\varphi$. Thus
$\psi(x,\bar c_0)\in\tp_K(e/D\cap K)$ divides over $B\cap K$.
This proves the claim.

It follows that, for every small tree-closed $B$ containing $b$,
\[
\SU_M(e/B)=\SU_K(e/B\cap K).
\tag{**}
\]
Here $\SU_K$ is computed in $T_\sigma$. For completeness, compare the
two ranks by induction on their finite lower bounds. In the ambient
structure, a parameter extension may be replaced by its tree closure
without changing the type's rank or forking, since tree closure is
contained in definable closure. Claim $(*)$ then sends every ambient
forking extension to a local forking extension. Conversely, a local
extension to a small set $F\subseteq K$ containing $B\cap K$ is
obtained in $M$ by passing to $\tcl(BF)$, and
$
\tcl(BF)\cap K=F.
$
Claim $(*)$ applies in this direction as well. These observations give
both steps of the induction. Since the local rank is finite, they
give $(**)$ and finiteness of the ambient rank.

Now let $A=\tcl(Ca)$. Since $a$ and an enumeration of $A$ are
interdefinable over $C$,
\[
\SU_M(a/C)=\SU_M(A/C).
\]
Choose a block enumeration
$
C=A_0\subseteq A_1\subseteq\cdots\subseteq A_m=A
$
as in Lemma~\ref{lem:block-enumeration}, with
$A_i=\tcl(A_{i-1}e_i)$ and
$\operatorname{pred}(e_i)\in A_{i-1}$. Let $K_i=\Comp(e_i)$.
The other points of $A_i\setminus A_{i-1}$ are forced
$1$-successors and are definable over $A_{i-1}e_i$. Hence $(**)$ gives
\[
\SU_M(A_i/A_{i-1})
=\SU_M(e_i/A_{i-1})
=\SU_{T_{\sigma(K_i)}}(e_i/A_{i-1}\cap K_i).
\]
All these ranks are finite. The Lascar inequalities therefore give
\[
\SU_M(A/C)
=\sum_{i=1}^m\SU_M(A_i/A_{i-1}).
\]
For a fixed component $K$, the elements $e_i$ with $K_i=K$ enumerate
$A_K\setminus C_K$, and $A_{i-1}\cap K$ consists of $C_K$ together
with the elements of this tuple already added. Applying finite-rank
additivity in $T_{\sigma(K)}$ and grouping the summands by component,
we obtain
\[
\begin{aligned}
\SU_M(a/C)
&=\sum_{i=1}^m
   \SU_{T_{\sigma(K_i)}}(e_i/A_{i-1}\cap K_i)\\
&=\sum_{K\in\mathcal K(A/C)}
   \SU_{T_{\sigma(K)}}(\bar a_K/C_K),
\end{aligned}
\]
as desired.

To compute the rank of $P_\tau$, take $C=\tcl(\emptyset)$. A point
of $P_\tau$ adds exactly one infinity element at each infinity-node
on the path to $\tau$. Thus the formula just proved gives the upper
bound
\[
\SU_M(P_\tau)
\leq\sum_{\substack{\sigma\in I(\Gamma)\\\sigma\leq\tau}}
          \SU(T_\sigma).
\]
For the reverse inequality, choose, successively along this path, an
element of maximal local rank in each component. Such elements exist
in the saturated components, and the definition of the canonical tree
product permits all these choices. The resulting point has rank
equal to the displayed sum. Since the plan is finite and
$M=\bigcup_{\tau\in\Gamma}P_\tau$, the formula for $\SU(M)$ follows.
In particular, $\Th(M)$ is supersimple of finite $\SU$-rank.

Finally, suppose that each $T_\sigma$ is one-based. Their finite
many-sorted disjoint union is one-based and of finite rank, and
adding the finitely many definable product sorts preserves these
properties. One-basedness is preserved under interpretations of
finite-rank simple theories, by the finite-rank reduct theorem of
N\"ubling.\footnote{H.~N\"ubling, \emph{Reducts of stable,
CM-trivial theories}, Journal of Symbolic Logic \textbf{70} (2005),
1025--1036, Section~3.}
Lemma~\ref{lem:coordinate-presentation} therefore implies that
$\Th(M)$ is one-based.

For an admissible NIC structure, apply
Theorem~\ref{thm:characterization}. The structure and its canonical
tree product are interdefinable, so the conclusions transfer.
\end{proof}

\begin{remark}
For canonical tree products, Theorem~\ref{thm:rank-one-based} does not
require trivial algebraic closure in the components. When their
algebraic closure is nontrivial, the formula still uses the local ranks
$\SU_{T_{\sigma(K)}}(\bar a_K/C_K).$
In the admissible case with trivial component algebraic closure, these
are the rank contributions of the new component coordinates.
\end{remark}

\subsection{Tree envelopes}

Finite homogeneous component substructures give finite approximations of
the tree product. We use these to construct envelopes for smooth
approximation.

\begin{definition}
Let $E$ be a finite substructure of a structure $N$.  We say that $E$ is
\emph{fully homogeneous in $N$} if the following two conditions hold.
\begin{enumerate}
\item Every $\emptyset$-definable relation of $N$ induces an
$\emptyset$-definable relation on $E$.
\item For all tuples $\bar a,\bar b$ from $E$ of the same length, with length at
most $|E|$,
$$
\tp_E(\bar a)=\tp_E(\bar b)
   \quad\Longleftrightarrow\quad
   \tp_N(\bar a)=\tp_N(\bar b).
$$
\end{enumerate}
\end{definition}

We use full homogeneity in this sense throughout the subsection.
Recall that a structure is smoothly approximable if it is
$\aleph_0$-categorical and each finite subset is contained in a finite
substructure fully homogeneous in the ambient structure.

\begin{definition}
Let
$
M=\prod_{\sigma\in I(\Gamma)}^{\!*}N_\sigma
$ 
be a canonical $\Gamma$-tree product.  For each $\sigma\in I(\Gamma)$, let $E_\sigma\subseteq N_\sigma$
be a finite substructure.  The corresponding \emph{tree envelope} $E_\Gamma((E_\sigma)_{\sigma\in I(\Gamma)})$
is the finite substructure of $M$ induced on the set of all tree nodes whose
$\sigma$-coordinate lies in $E_\sigma$ whenever that coordinate is defined.
Equivalently, for each $\tau\in\Gamma$,
$ 
P_\tau(E_\Gamma)
   =
   \prod_{\substack{\sigma\in I(\Gamma)\\ \sigma\leq\tau}} E_\sigma,
$ 
with the empty product interpreted as a singleton.
\end{definition}

\begin{theorem}
\label{thm:tree-envelopes}
Let
$ 
M=\prod_{\sigma\in I(\Gamma)}^{\!*}N_\sigma
$  
be a canonical $\Gamma$-tree product.  Suppose that, for each
$\sigma\in I(\Gamma)$, the component structure $N_\sigma$ has a cofinal family of
finite fully homogeneous substructures.

Then the finite tree envelopes $E_\Gamma((E_\sigma)_{\sigma\in I(\Gamma)})$
are cofinal among finite subsets of $M$.  Moreover, if each $E_\sigma$ is fully
homogeneous in $N_\sigma$, then the tree envelope $E_\Gamma((E_\sigma)_{\sigma\in I(\Gamma)})$
is fully homogeneous in $M$.

Consequently, if each component structure $N_\sigma$ is smoothly approximable,
then the canonical $\Gamma$-tree product $M$ is smoothly approximable.
\end{theorem}

\begin{proof}
First let $F\subseteq M$ be finite.  For each $\sigma\in I(\Gamma)$, let $F_\sigma\subseteq N_\sigma$
be the finite set of all $\sigma$-coordinates appearing in elements of $F$.  By
cofinality, choose a finite fully homogeneous substructure $E_\sigma\subseteq N_\sigma$
containing $F_\sigma$.  The resulting tree envelope contains $F$. This proves cofinality.

Now suppose that each $E_\sigma$ is fully homogeneous in $N_\sigma$.
Let $E=E_\Gamma((E_\sigma)_{\sigma\in I(\Gamma)}).$
We show that $E$ is fully homogeneous in $M$.

Let $\bar a,\bar b$ be tuples from $E$ of the same length. In both $E$
and $M$, the coordinate presentation describes their types by the pure
tree diagram and the local types in the active components. The pure
tree diagram is unchanged. In a component $K$ of plan-type $\sigma$,
the local tuples lie in $E_\sigma$. By full homogeneity of $E_\sigma$
in $N_\sigma$, equality of their types in $E_\sigma$ is equivalent
to equality of their types in $N_\sigma$.

Therefore, for all tuples $\bar a,\bar b$ from $E$,
$$
\tp_E(\bar a)=\tp_E(\bar b)
   \quad\Longleftrightarrow\quad
   \tp_M(\bar a)=\tp_M(\bar b).
$$
This proves the tuple-homogeneity part.

We now verify the first clause of full homogeneity. Let $R$ be an $\emptyset$-definable $r$-ary relation of $M$. If two tuples from $E$
have the same type in $E$, they have the same type in $M$ by what we
have just proved. Hence they either both satisfy $R$ or both fail to
satisfy it. Thus $R\cap E^r$
is a union of complete $r$-types in $E$, or equivalently a union of
$\Aut(E)$-orbits. Since $E$ is finite, this union is
$\emptyset$-definable in $E$. This proves the required clause.

Thus $E$ is fully homogeneous in $M$.

Suppose now that each $N_\sigma$ is smoothly approximable. The
cofinality argument gives a finite fully homogeneous tree envelope
containing any given finite $F\subseteq M$. Also, each component is
$\aleph_0$-categorical. By Lemma~\ref{lem:coordinate-presentation}
and finiteness of the plan, $M$ is $\aleph_0$-categorical.
Therefore $M$ is smoothly approximable.
\end{proof}

\begin{remark}
These tree envelopes are obtained from finite homogeneous component
approximations. They give the nil-interaction version of the envelope
construction, and should be distinguished from Cherlin--Hrushovski
$\mu$-envelopes for general Lie coordinatization. Here the construction
is simpler because different components are orthogonal.
\end{remark}

\section{Finite tree products and multidimensional asymptotic classes}
\label{sec:mac-trees}

We work in the multidimensional framework of
\cite{AnscombeMacphersonSteinhornWolf2024}, using the connection between
exact classes and smooth approximation in \cite{Wolf2020}. The arguments
in this section are independent of Theorem~\ref{thm:main}.

Roughly speaking, a multidimensional asymptotic class has finitely many
approximate counting functions for each definable family, and the
parameters giving each function form a definable set. For a
multidimensional exact class, these functions give the cardinalities
exactly.

\subsection{Finite tree products}

We allow a different finite component structure at each infinity-node.
For each $\sigma\in I(\Gamma)$, let $A_\sigma$ be a finite
$L_\sigma$-structure, with the languages $L_\sigma$ pairwise disjoint.
Fix a nonempty coordinate sort $U_\sigma(A_\sigma)$
of $A_\sigma$.  In the one-sorted case, this is just the universe of $A_\sigma$.

Let
$\bar A=(A_\sigma)_{\sigma\in I(\Gamma)}.$
The finite tree product $\Gamma(\bar A)$ is defined as in
Section~\ref{sec:trees}, except that the coordinate chosen at an infinity-node
$\sigma$ is taken from the coordinate sort $U_\sigma(A_\sigma)$.  Thus
$$
|P_\tau(\Gamma(\bar A))|
   =
   \prod_{\substack{\sigma\in I(\Gamma)\\ \sigma\leq\tau}}
   |U_\sigma(A_\sigma)|.
$$
The relations on the coordinate sort of $A_\sigma$ are interpreted
separately in each copy of that component. If a component language has
function symbols, we replace them by their graphs. This changes the
notation but not the definable sets.

We will also use the enriched many-sorted structure
$\Gamma^{\mathrm{en}}(\bar A).$
This structure has the tree-product sort and the component sorts of
each $A_\sigma$. It also has coordinate maps from $P_\tau$ to
$U_\sigma(A_\sigma)$ for each infinity-node $\sigma$ on the path to
$\tau$. The structures $A_\sigma$ are named on their sorts, and the
local relations on tree nodes are pulled back along these maps.

\begin{proposition}\label{prop:pure-tree-mec}
Fix a finite tree plan $\Gamma$.  Let $\mathcal T^{\mathrm{en}}_\Gamma$ be the
class of enriched finite pure tree products obtained by assigning a finite nonempty
pure set $X_\sigma$ to each $\sigma\in I(\Gamma)$.  Then
$\mathcal T^{\mathrm{en}}_\Gamma$ is a polynomial multidimensional exact class.
The counting functions are integer-valued polynomial functions in the variables
$
N_\sigma=|X_\sigma|
   \,\, \text { for }\sigma\in I(\Gamma).
$
Its one-sorted tree-product reduct is a weak polynomial multidimensional exact
class.  If the coordinate enrichment is uniformly $\emptyset$-definable, or more
generally uniformly $\emptyset$-interpretable so that the reduct and enrichment are
uniformly bi-interpretable, then the word ``weak'' can be omitted.
\end{proposition}

\begin{proof}
First consider the many-sorted disjoint union of the pure sets
$X_\sigma$. It forms a polynomial multidimensional exact class.
Indeed, formulas are Boolean combinations of equalities, and we
partition the parameter space according to equality patterns in each
sort. On each part, the cardinality is an integer-valued polynomial
in the variables
$N_\sigma=|X_\sigma|.$

We define the enriched tree product uniformly in this disjoint union.
For each $\tau$, take the product of the sorts $X_\sigma$ with
$\sigma\in I(\Gamma)$ and $\sigma\leq\tau$. If there is no such
$\sigma$, take a singleton. The tree maps are given by projections
and extensions of coordinates, and the coordinate maps are projections.
This construction uses no quotients.

Conversely, each $X_\sigma$ is a named sort of the enrichment.
Hence the enriched class is uniformly bi-interpretable with a
definitional expansion of the disjoint union. The result follows
from preservation of multidimensional exact classes under uniform
bi-interpretability.

For the one-sorted reduct, discard the component sorts and coordinate
maps. This gives a uniform interpretation, so the reduct is a weak
multidimensional exact class. If the enrichment is uniformly recoverable,
we have a uniform bi-interpretation and retain the definability clauses.
\end{proof}

\begin{theorem}
\label{thm:finite-tree-mac}
Fix a finite tree plan $\Gamma$.  For each $\sigma\in I(\Gamma)$ let
$\mathcal C_\sigma$ be an $R_\sigma$-multidimensional asymptotic class of finite
$L_\sigma$-structures, with pairwise disjoint languages and with a distinguished
nonempty coordinate sort $U_\sigma$.  Let $\mathcal C^{\mathrm{en}}_\Gamma$ be
the class of enriched finite tree products
$
\Gamma^{\mathrm{en}}\bigl((A_\sigma)_{\sigma\in I(\Gamma)}\bigr),
   \,\, A_\sigma\in\mathcal C_\sigma.
$
Then $\mathcal C^{\mathrm{en}}_\Gamma$ is an
$R_\Gamma$-multidimensional asymptotic class, where $R_\Gamma$ is generated,
under addition and multiplication, from the functions in the $R_\sigma$ and the
coordinate-size functions
$N_\sigma(A_\sigma)=|U_\sigma(A_\sigma)|.$
If every $\mathcal C_\sigma$ is a multidimensional exact class, then
$\mathcal C^{\mathrm{en}}_\Gamma$ is a multidimensional exact class.

The one-sorted tree-product reduct is uniformly interpretable in
$\mathcal C^{\mathrm{en}}_\Gamma$, and hence is a weak multidimensional
asymptotic, respectively exact, class.  If the coordinate enrichment is uniformly
$\emptyset$-definable, or more generally uniformly $\emptyset$-interpretable so
that the reduct and enrichment are uniformly bi-interpretable, then the word
``weak'' can be omitted.
\end{theorem}

\begin{proof}
We first form the disjoint many-sorted union of the classes
$\mathcal C_\sigma$. Multidimensional asymptotic and exact classes
are closed under finite unions of this form. The resulting measuring
functions are obtained from the $R_\sigma$ by addition and multiplication.

We interpret the enriched tree product uniformly in this disjoint union.
For $P_\tau$, take the product of the sorts $U_\sigma(A_\sigma)$
along the path to $\tau$. The tree structure and coordinate maps are
defined using projections and extensions of coordinates. We keep the
component relations on the named sorts and pull them back to the
appropriate tree nodes.

Since the component sorts are named in the enrichment, the converse
interpretation is immediate. Thus we have a uniform bi-interpretation
with a definitional expansion of the disjoint union. Preservation under
uniform bi-interpretability proves the first assertion.

We now check the measuring functions. The interpretation uses no
quotients: each new sort is a finite disjoint union of finite products
of coordinate sorts. Therefore addition and multiplication of the
component measuring functions and coordinate-size functions suffice.
Using only the general interpretability theorem would give the weaker
conclusion with $\operatorname{Frac}(R_\Gamma)$ in place of $R_\Gamma$.

Discarding the component sorts and coordinate maps gives a uniform
interpretation of the one-sorted reduct. The interpretability theorem
therefore gives the weak conclusion. If the enrichment can be recovered
uniformly, the interpretation is a bi-interpretation, and the definability
clauses follow as well.
\end{proof}

\begin{corollary}
\label{cor:generalized-measurable-tree-products}
Let $\mathcal C^{\mathrm{en}}_\Gamma$ be as in
Theorem~\ref{thm:finite-tree-mac}.  Let
$(M_i)_{i\in J}$
be a family of members of $\mathcal C^{\mathrm{en}}_\Gamma$, and let $\mathcal U$
be an ultrafilter on $J$ such that the cardinalities $|M_i|$ tend to infinity modulo
$\mathcal U$.  Then the ultraproduct
$\prod_{i\in J} M_i/\mathcal U$
is generalized measurable in the sense of
\cite{AnscombeMacphersonSteinhornWolf2024}.  If the component classes are
multidimensional exact classes, then the ultraproduct is ring-measurable.
\end{corollary}

\begin{proof}
Apply the ultraproduct theorem for multidimensional asymptotic classes
and its exact-class version.
\end{proof}

\begin{example}
Let $\mathcal C_\sigma$ consist of the two-sorted structures $(F,V)$,
where $F$ is a finite field and $V$ is a finite-dimensional vector space
over $F$, in the usual language with scalar multiplication. Use $V$ as
the tree-coordinate sort and retain $F$ in the enrichment. These
vector-space classes are polynomial multidimensional asymptotic classes.
By Theorem~\ref{thm:finite-tree-mac}, their enriched finite tree products
form a multidimensional asymptotic class. The corresponding ultraproducts
are generalized measurable.

Local algebraic closure here is linear span. Hence this example is
outside the admissible NIC theorem with trivial component closure.
It suggests considering a closure which combines tree closure with
algebraic closure in the components.
\end{example}

\begin{remark}
The independent component sizes explain the use of multidimensional
asymptotics in Theorem~\ref{thm:finite-tree-mac}. Even for pure trees,
the parameters $N_\sigma=|U_\sigma(A_\sigma)|$
may grow at different rates. The multidimensional framework allows
these rates to remain independent, whereas ordinary asymptotic classes
use a single growth scale.
\end{remark}

\section{Further directions}

The notion of coordinatization used here is strong. Different components
have no interaction beyond that supplied by the tree. This assumption
gives the normal form for finite types and the formulas for dimension
and measure. It is natural to ask how far it can be weakened.

First, one can allow some relations between different components. The
generic bipartite graph shows what must be addressed: a parameter in
one part can change the type of a point in the other part without
changing its local type. Thus a more general measuring construction
would need information about these additional relations. One possible
direction is to study their connection with higher amalgamation in
pseudofinite $\aleph_0$-categorical theories; see \cite{Kruckman2019}.

Second, one can allow nontrivial algebraic closure in the components.
Vector spaces and measurable modules are natural examples to consider.
For such structures, we would replace $\tcl$ by a closure combining
tree closure and local algebraic closure. The measuring construction
would then use local fibers over algebraic bases in place of the
independent singleton choices used here.

\section*{Acknowledgements}

This paper is a revised version of the author's 2021 preprint.
The author would like to thank Camryn Hill for her guidance
and support, and Alex Kruckman and Philip Scowcroft for
their helpful comments.

\begin{samepage}
\medskip\noindent{\textbf{AI Use  Statement:}} ChatGPT was used to help improve the language and presentation of the manuscript, and to assist with proof checking and revision. The author takes full responsibility for the content of the manuscript.
\end{samepage}

\end{document}